\documentclass[12pt,leqno]{amsart}
\usepackage{amsmath,amsthm}
\usepackage{amsfonts}
\newcommand{\cal}{\mathcal}

 \newtheorem{res}{Result}[section]
 \newenvironment{remark}{\noindent\stepcounter{res}\textbf{Remark \theres}}{\medskip}
 \newtheorem{theorem}[res]{Theorem}
 \newtheorem{proposition}[res]{Proposition}
 \newtheorem{lemma}[res]{Lemma}
 
 \newtheorem{definition}[res]{Definition}
\newtheorem{conjecture}[res]{Problem}
\numberwithin{equation}{section}
\newcommand{\half}{{\mbox{\small $\frac{1}{2}$}}}
\def\Kappa{{\mathcal K}}
\begin{document}
\bibliographystyle{plain}
\title{RANDOM ORDERINGS OF THE INTEGERS AND CARD SHUFFLING}
\author{Saul Jacka}
\address{Department of Statistics, University of Warwick, Coventry CV4 7AL, UK}
\email{s.d.jacka@warwick.ac.uk}
\author{Jon Warren}
\email{j.warren@warwick.ac.uk}
\thanks{{\it 2001 MSC Subject classifications:} Primary: 60B15; Secondary: 60G09, 37A40, 37H99, 60J05, 03E10}
\thanks{{\it Keywords and phrases:} RIFFLE SHUFFLES, QUASI-UNIFORM MEASURES, EXCHANGEABLE ORDERINGS, SHUFFLE IMBEDDING SHUFFLE}
\thanks{Abbreviated Title: RANDOM ORDERINGS AND CARD SHUFFLING}

\begin{abstract}
In this paper we study  random
orderings of the integers with a certain invariance property.
We
 describe all such orders in a simple way. We define and represent random shuffles of
a countable set of labels
 and then give an 
 interpretation of these orders 
in terms of a class of generalized riffle shuffles.
\end{abstract}
\maketitle

\section{Introduction}
In Jacka and Warren (1999) we defined deterministic shuffles on (a countable set of
labels indexed by) $\mathbb N$. In this paper we define random shuffles
on $\mathbb N$ and represent their laws in terms of the laws of pairs of
random variables with uniform marginals (Theorem \ref{represent}).
A natural subclass of random shuffles are the shuffle imbedding
shuffles:
those shuffles whose restrictions to \(\{1,\ldots,n\}\) induce a random
walk on ${\mathcal S}_n$).
Partly in order to study such shuffles, and partly because they are of
substantial interest in their own right, we introduce and study the class of 
$\cal I$-invariant orderings: random orderings of $\mathbb Z$ whose laws
are invariant under increasing relabellings. Section 3 is devoted to defining and representing
$\cal
I$-invariant orderings in terms of {\it quasi-uniform} measures (Theorem
\ref{orderth}).

\section{Preliminaries}

We denote by $\cal O$ the class of all strict total orderings
of $\mathbb Z$. This inherits a natural measurable structure
as a subset of $2^{{\mathbb Z}\times {\mathbb Z}}$. We will
denote a generic element of $\cal O$ by $ \lhd$, and write $m
\lhd n $ if $m$ is less than $n$ under $\lhd$.

Given any strictly increasing map $f: {\mathbb Z} \mapsto
{\mathbb Z}$ there is a naturally induced map $\hat{f}:{\cal
O} \mapsto {\cal O}$ defined by \[ m \stackrel{\hat{f}}{\lhd}
n   \qquad \text{ if and only if } \qquad f(m) \lhd f(n), \]
where we are denoting the image of the ordering $\lhd$ under
$\hat{f}$ by $\stackrel{\hat{f}}{\lhd}$ \begin{definition}

A probability measure ${\mathbb P}$ on  $\cal O$ is said to be
$\cal I$-invariant if $ {\mathbb P} \circ
\hat{f}^{-1}={\mathbb P}$, for all strictly increasing $f$.
\end{definition}

Our purpose is to give an explicit description of all such
invariant random order relations.

Hirth and Ressel (2000) considered
a similar invariance property on
 random orderings where $f$ ranges over finite permutations.
 Such random orderings are called exchangeable. Their
 characterization of the laws of these orderings  is very
 reminiscent of our results on ${\cal I}$-invariant orderings.
 This is related to the well-known fact (see Lemma
 \ref{exlemma2} below) that exchangeability of an infinite
 sequence of random variables is equivalent to an apparently
 weaker condition involving the action of increasing maps.

The following is a fundamental example that illustrates the
connection with riffle shuffles.
 Suppose that $\bigl( Z_i\bigr)_{i \in {\mathbb Z}}$
is a doubly infinite sequence of independent, identically
distributed random variables taking values in $\bigl\{
0,1\bigr\}$. Define a random order $\lhd$ as follows.
\begin{equation*} m \lhd n \qquad \text{ if and only if 
either} \qquad (Z_m=Z_n\quad \text{and}\quad m<n) \quad \text{ or } \quad Z_m<Z_n.
\end{equation*} In effect we split $\mathbb Z$ into two
equivalence classes, namely $\bigl\{n: Z_n=0\bigr\}$ and
$\bigl\{n:Z_n=1\bigr\}$; we order the former below the latter,
and within each class we preserve the natural order. This is
an example of an $\cal I$-invariant ordering.

Perhaps the most celebrated example of a shuffle is the
Gilbert-Shannon-Reeds riffle shuffle. As described in Diaconis
(1998) this is obtained by the following  recipe.
Suppose that  our cards are labelled $1,2 \ldots n$ and  we
take $n$ independent random variables $U_1, \ldots U_n$ each
uniformly distributed on $[0,1]$. Order the cards initially so
that card $k$ is above card $l$ whenever $U_k>U_l$.  Then 
reorder the cards according to the values of $2U_k \mod 1$.
The random permutation that must be applied to reorder the
cards is the GSR shuffle. When we look at the inverse
permutation  we see that the cards are divided into two
subpacks according to whether $U_k$ is smaller or greater than
one-half, and within each subpack order is preserved. Thus the
inverse of the GSR shuffle is described by (the restriction to
$\{1, \ldots n\}$ of) an $\cal I$-invariant ordering. In Bayer
and Diaconis (1992), this description is used to
investigate the speed of mixing of the GSR shuffle.
 Generalisations based  on replacing $u \mapsto 2u
\mod 1$ with other  maps have been studied by other authors,
see for example Lalley (1999). We will see that there
is a quite general correspondence to be made between  $\cal
I$-invariant orderings, families of riffle shuffles and a
class of measure-preserving maps on $[0,1]$.

\section{Describing $\cal I$-invariant orderings}

\begin{definition}
A probability measure $\mu$ on $[0,1]$ is called quasi-uniform
if it satisfies \[ \mu\bigl\{x \in [0,1]: \mu[0,x)\leq x \leq
\mu[0,x] \bigr\}=1. \] \end{definition}

\begin{remark}
It is not hard to show that the set of all quasi-uniform measures is
closed with respect to the topology of waek convergence of probability
measures on $[0,1]$.
\end{remark}

Such a measure is really quite a simple object, and it may be
described as follows.  \begin{lemma} \label{dlemma} Suppose
that $F$ is a closed subset of $[0,1]$, and $\lambda_F$ is the
measure with density $1_F$ with respect to Lebesgue measure.
Corresponding to each open component $G_i$ of its complement
$F^c$ is a point mass $m_i \delta_{x_i}$, of size $m_i$ equal
to the length of the interval $G_i$, situated at position
$x_i$, which is either the left or right hand end of $G_i$.
Then the measure $\mu$ given by \[ \mu= \lambda_F + \sum_i m_i
\delta_{x_i}, \] is quasi-uniform. Moreover every
quasi-uniform $\mu$ can be decomposed in this fashion, and the
decomposition is unique (up to the labelling of the
intervals). \end{lemma} We omit the proof which is elementary.

Notice that it is possible for two distinct masses in the
above decomposition to be placed at the same point.

Suppose $\mu$ is quasi-uniform then the measure
 $\mu^\prime$,  obtained by
inverting the distribution function  of $\mu$:
\[
\mu^\prime[0,y]= \inf \{x: \mu[0,x] \geq y\},
\]
is also quasi-uniform.
This corresponds to switching each mass $m_i$ in the 
decomposition of $\mu$ to being at the opposite end of the
interval to which it belongs. If $X$ and $Y$ are random
variables on the same probability space with the law of $X$
being $\mu$ and the law of $Y$ being $\mu^\prime$ and so that
$X$ and $Y$ are equal or take values at either end of a
component of $F^c$ then let us say that such $X$ and $Y$ form
a conjugate pair. Notice that their  joint law  is specified 
completely by the above description. Such a pair may contain a
little more information than either variable separately:
whenever they are not equal  they together determine an
interval of $F^c$. 

Now for any quasi-uniform $\mu$ we construct an ${\cal
I}$-invariant ordering whose law we denote by ${\mathbb
P}^{\mu}$. Consider an infinite sequence of independent pairs
of random variables $\bigl(X_n,Y_n\bigr)_{n \in {\mathbb Z}}$.
For each $n$ the variables $X_n$ and $Y_n$  form a conjugate
pair with $X_n$ distributed according to $\mu$. Now, supposing
that $m<n$ (with respect to the natural ordering of the
integers) take $m \lhd n$ if and only if  one of the following
happens: \begin{equation} \label{condition} \begin{cases}
X_m<X_n,&\\ Y_m\kern2pt <\kern2pt Y_n,&\\ X_m=X_n\phantom{,}&\kern-10pt >Y_m=Y_n. \end{cases}
\end{equation}

It is easy to check that this works and defines an $\cal
I$-invariant ordering. Moreover the strong law of large
numbers implies that \begin{equation} \label{limits}
\begin{cases} X_n= &\lim\limits_{N \rightarrow \infty} \frac{1}{N}
\sum_{-N}\limits^n 1_{(k\lhd n)}
 \\
Y_n= &\lim\limits_{N \rightarrow \infty} \frac{1}{N} \sum\limits_n^N
1_{(k\lhd n)} , \end{cases} \end{equation} exist almost
surely. Notice that, since $\mu$ can be recovered from $\lhd$
as the empirical distribution of the sequence $X_n$, ${\mathbb
P}^{\mu_1} \neq {\mathbb P}^{\mu_2}$ if $\mu_1 \neq \mu_2$. The
following theorem says that by taking mixtures of  orderings
of this form we obtain all possible $\cal I$-invariant
orderings. \begin{theorem} \label{orderth} Suppose that $\lhd$
is an $\cal I$-invariant ordering. Then almost surely, the
random variables defined by \eqref{limits} exist, and for any
$m$ and $n$ the relation $m\lhd n$ holds if and only if
\eqref{condition} does.

Moreover  the sequence  of random variables

$\bigl(X_n\bigr)_{n \in{\mathbb Z}}$  is exchangeable and with
probability one  it admits  an  empirical distribution
$\mu(X)$ which is quasi-uniform.

Conditional on
$\mu(X)=\mu$ the law of $\lhd$ is ${\mathbb P}^\mu$. 
\end{theorem}

In general any  ordering $\lhd$ belonging to $\cal O$ 
projects to an equivalence relation, $\sim$, on $\mathbb Z$  defined
by
\begin{displaymath}
n \sim m \Leftrightarrow \hbox{ there are only finitely many $k$
between (with respect to $\lhd$) $n$ and $m$.}\end{displaymath}
If the ordering
is $\cal I$-invariant then it follows from the above theorem
that this partition is exchangeable in the sense studied by  Kingman
(1982).

The proof of Theorem \ref{orderth} hinges on  the elementary
observation of the next lemma, which begins to explain  the
role of quasi-uniform measures. \begin{lemma}
\label{quasilemma} Suppose that $\lhd$ is some fixed ordering.
Define, for each $n$, \begin{eqnarray*} \overline{X}_n &=
\limsup\limits_{N \rightarrow \infty} \frac{1}{N} \sum\limits_1^N 1_{(k\lhd
n)} \\ \underline{X}_n &= \liminf\limits_{N \rightarrow \infty}
\frac{1}{N} \sum\limits_1^N 1_{(k\lhd n)} \end{eqnarray*} Suppose
that the measures $\nu^{(N)}$, defined by \[
\nu^{N}[0,x]=\frac{1}{N} \sum_ 1^N 1_{(\overline{X}_k \leq x)}
\qquad x\in[0,1], \] converge weakly to a probability measure
$\nu$ as $N$ tends to infinity. Then \[ \nu[0,\overline{X}_1)
\leq \underline{X}_1 \leq \overline{X}_1 \leq
\nu[0,\overline{X}_1]. \] \end{lemma} \begin{proof} It is  an
easy consequence of the transitivity of $\lhd$ that: \[
\overline{X}_k <\overline{X}_1 \Rightarrow k \lhd 1
\Rightarrow \overline{X}_k \leq \overline{X}_1, \] for any $k
\in {\mathbb Z}$. Thus \[ \frac{1}{N}\sum_ 1^N
1_{(\overline{X}_k < \overline {X}_1)} \leq \frac{1}{N}
\sum_1^N 1_{(k\lhd n)} \leq \frac{1}{N} \sum_ 1^N
1_{(\overline{X}_k \leq \overline {X}_1)} . \] But the
left-hand side is $\nu^{(N)}[0,\overline{X}_1)$ while the
right-hand side is $\nu^{(N)}[0,\overline{X}_1]$ and by virtue
of weak convergence: \begin{eqnarray*}
\nu[0,\overline{X}_1)&\leq \liminf
\nu^{(N)}[0,\overline{X}_1)\\ \nu[0,\overline{X}_1]& \geq
\limsup \nu^{(N)}[0,\overline{X}_1]. \end{eqnarray*}
\end{proof} \begin{lemma} \label{exlemma1} Suppose that $\lhd$
is an ${\cal I}$-invariant ordering. Then the family of random
variables \[ \bigl( 1_{(k \lhd 0)}; k >0 \bigr) \] is
exchangeable. \end{lemma} \begin{proof} It suffices to check
that for finite collections of positive integers $j_1 \ldots
j_m$ and $k_1 \ldots k_n$ the value of \[ {\mathbb E}\left[
1_{( 0 \lhd j_1)} \ldots 1_{( 0 \lhd j_m)} 1_{( k_1 \lhd 0)}
\ldots 1_{( k_n \lhd 0)} \right] \] depends only on $n$ and
$m$. Now  replace $ 1_{( k_i \lhd 0)}$ by $1-1_{( 0 \lhd
k_i)}$, multiply out and apply ${\cal I}$-invariance to obtain
an expression involving terms: ${\mathbb
  E}[1_{(0 \lhd 1)}1_{(0\lhd 2)} \ldots 1_{(0 \lhd k)}]$ for
  $m \leq k
\leq m+n$.
\end{proof}
This lemma is actually a special case of  the next result, for
the proof of  which we refer the  reader to Aldous(1985).
\begin{lemma} \label{exlemma2} Suppose that a sequence of
random variables $\bigl(X_k ;k\in {\mathbb Z} \bigr)$ is
${\cal I}$-invariant, in the sense that for any increasing
function $f:{\mathbb Z} \mapsto {\mathbb Z} $ \[ \bigl(X_k
;k\in {\mathbb Z} \bigr)\stackrel{law}{=}\bigl(X_{f(k)}
;k\in{\mathbb Z} \bigr) \] then in fact  $\bigl(X_k ;k\in
{\mathbb Z} \bigr)$  are exchangeable- the sequence admits
with probability one an empirical distribution and conditional
on it the random variables are independent and identically
distributed. \end{lemma}

\begin{proof}[Proof of theorem \ref{orderth}]

We begin by observing that the variables  $(X_k,Y_k)$ exist by
virtue of  the exchangeability property of  Lemma
\ref{exlemma1} and De Finetti's Theorem. Moreover the law of
the sequence  of pairs $(X_k, Y_k)$ is ${\cal I}$-invariant so
we may deduce from Lemma \ref{exlemma2} that it is, in fact, an
exchangeable sequence. It follows from Lemma \ref{quasilemma}
that the  empirical distributions for both $X_k$ and $Y_k$
must be quasi-invariant.

The next step is to show that the variables $(X_k,Y_k)$
determine the ordering $\lhd$. Divide ${\mathbb Z} /\{0\}$
into three classes. \begin{align*} U_0 &= \{k: \text{ either }
X_k > X_0 \text{ or } Y_k>Y_0\} \\ E_0 &= \{k:  X_k = X_0
\text{ and } Y_k=Y_0á\} \\ B_0 &= \{k: \text{ either } X_k <
X_0 \text{ or } Y_k<Y_0\}. \end{align*} Notice that, since $j
\lhd k$ implies $X_k \geq X_j$ and $Y_k \geq Y_j$, we must have
any element of $U_0$ ordered above any element of $E_0$ which
in turn must be ordered above any element of $B_0$. Because of
the exchangeability of $X$ and $Y$, the three classes have
limiting sizes: \[ \vert U_0 \vert= \lim \frac{1}{N} \sum_1^N
1_{(k \in U_0)}=\lim \frac{1}{N} \sum_{-N}^{-1} 1_{(k \in
U_0)} \] and similarly for $\vert E_0 \vert$ and $\vert B_0
\vert$. The exchangeability of $X$ and $Y$ also implies that
if the size $\vert E_0 \vert$ of $E_0$ is zero then it is, in 
fact, empty. Otherwise the empirical distributions of $X$ and
$Y$ have atoms at the values of $X_0$ and $Y_0$. We claim  the
restriction of $\lhd $ to $E_0$ either preserves or reverses
the natural order: it then follows that in the former case:
$Y_0=X_0-\vert E_0 \vert$,  while in the latter case: 
$Y_0=X_0+\vert E_0 \vert$. This  then  establishes that the
ordering $\lhd$ is determined by the sequence $(X_k,Y_k)$
according to \eqref{condition}. 

To prove the claim  of the previous paragraph suppose that
$0<j<k$, and let $p$ be the  probability  ${\mathbb P}(0 \lhd
k \lhd j, \text{ and } j \in E_0 )$. Now \[ \frac{1}{N-j}
\sum_{r=j+1}^N 1_{(0 \lhd r \lhd j, \text{ and } j \in E_0 )}=
\frac{1_{(0 \lhd j\text{ and } j \in E_0
  )}}{N-j} \sum_{r=j+1}^N \left( 1_{(r \lhd j  )}-1_{(r \lhd 0
   )}\right),
\]
must converge to $1_{(0 \lhd j \text{ and } j \in E_0
)}(Y_j-Y_0)=0$ in $L^1$ (by bounded convergence),
 yet its expectation is , for
all $N$, equal to $p$.  This, and similar versions show that
if $r \in E_0$ then the only $s$ between $0$ and $r$ with
respect to $\lhd$ are also between $0$ and $r$ in the natural
ordering. But now we may replace $0$ by $t$ in this statement,
then  by noting that if
 $t \in E_0$ then $E_t=E_0$ we deduce that whenever $r$ and
$t$  both belong to $E_0$ then $s$ being between them with
respect  to $\lhd$ implies $s$ is between them with respect to
the natural order. 

To complete the proof of the theorem condition on the joint
empirical measure of $(X,Y)$ to reduce to the iid case. The
arguments in the previous step show that each $X_k$ and $Y_k$
must form a conjugate pair, and that the conditional 
distribution of $\lhd$ is ${\mathbb P}^\mu$ where $\mu$ is the
empirical measure of $X$. 
 \end{proof}

Let us close this section by noting another natural invariance
property that one might impose on a random ordering: a
probability measure ${\mathbb P}$ on  $\cal O$ is said to be
$\cal T$-invariant if $ {\mathbb P} \circ
\hat{f}^{-1}={\mathbb P}$, for all $f$ of the form $f(n)=n+a$
for some $a \in {\mathbb Z}$.

\begin{conjecture} Obtain an explicit description of all $\cal
T$-invariant orderings. \end{conjecture}

${\cal T}$-invariant orders have a much richer structure than
${\cal I}$-invariant ones, as the following example
illustrates.
 Let $I_n$ be a stationary sequence of
$\{0,1\}$ random variables and construct $\lhd$ as follows.
\begin{itemize} \item If $I_n=I_m$ then $\lhd$ agrees with the
natural order. \item the upper class $\{I_n=1\}$ has slipped
one place relative to $\{I_n=0\}$. Thus for example if $n<m$
belong to class $0$ and class $1$ respectively then $n\lhd m$
unless there is no $n<k<m$ with $I_k=1$- if this happens then
$m \lhd n$. \end{itemize} Notice that making the random
variables $I_n$ independent does not make $\lhd$  ${\cal
I}$-invariant. Something more interesting is happening!

\section{Shuffling an infinite set of cards}

The first half of this section is based on  the (more
leisurely) account   contained in Jacka and Warren (1999)  of what it might
mean to shuffle an infinite set of cards. The state of an
infinite pack of cards will be represented by an ordering of
the natural numbers. The second half of the section considers
classes of Markov processes (indexed by discrete time)  taking
values in the space of such orderings and shows how one such
class is naturally associated with the class of ${\cal
  I}$-invariant orderings we have studied in the previous
  section.

Recall the standard model for shuffling cards: $n$ cards
carrying {\em labels} $1$ through to $n$ each have a distinct
{\em position} $1$ through to $n$ in the pack. We associate
the state of the pack with a permutation $\rho $ belonging to
the permutation group on $n$ objects ${\cal S}_n$. If
$\rho(k)=m$ then we say that the card carrying label $k$ is in
position $m$ in the pack.  A completely randomized pack simply
means choosing $ \rho$ according to the uniform measure on
${\cal S}_n$. A {\em shuffle} $S$ is a possibly random
permutation (belonging to ${\cal S}_n$!) of the positions in
the pack. Thus $S(m)= m^\prime$ means that the card that was
in position $m$ is moved to position $m^\prime$. Consequently
the state of the pack is changed from $\rho$ to $S\rho$. In
this way $S$ induces a map $\hat{S}: {\cal S}_n \mapsto {\cal
S}_n$ defined by $\hat{S}(\rho)= S\rho$. Such an $\hat{S}$
ignores the labelling of the pack. If $r$ ({\it also} belonging to 
${\cal S}_n$) is used to change the labels so that the card
that now carries the label $k$ is the card that previously
carried the label $r(k)$ and we denote by $\hat{r}$ the
induced map 
 $\hat{r}(\rho)=\rho r$ then we obtain the
commutation  relation \begin{equation} \label{com} \hat{S}
\circ \hat {r}= \hat{r} \circ \hat{S}, \end{equation} for all
$r \in {\cal S}_n$. Moreover any map $\hat{S}$  that commutes
with all relabellings is 
 induced by some  $S \in {\cal S}_n$.

We have rather laboured the point in the previous paragraph so
as to motivate our model for shuffling an infinite pack of
cards. We have seen that for a finite pack  the permutation
group plays three distinct roles - it describes the state of
the pack, it gives rise to shuffles, and it can be used to
relabel the pack. We proceed to the description of three
different objects that play these  roles in the infinite
framework. First note that for a finite pack we may also
specify the state of the pack by giving an ordering
$\lhd^{(n)}$ of $\{1, \ldots n\}$  related to our previous
description by means of a permutation $\rho$ via
\begin{equation} \label{correspond} k \lhd^{(n)} k^\prime \qquad
\text{ iff }\qquad \rho(k) <\rho(k^\prime). \end{equation}
With this approach we note that we may restrict the ordering
$\lhd^{(n)}$ to the first $n-1$ cards to obtain an ordering
$\lhd^{(n-1)}$. Moreover, if $\lhd^{(n)}$ is choosen uniformly
then $\lhd^{(n-1)}$ is uniformly distributed also. Because of
this consistency there is a unique measure $\lambda$ on the
space of total orderings of ${\mathbb N}$ so that the
restriction $\lhd^{(n)}$ of the ordering to $\{1, \ldots n\}$
is uniform. It is well known how to construct a random
ordering  distributed according to $\lambda$. Let $U_1,
U_2,\ldots $ be  an infinite sequence of independent random
variables uniformly distributed on $[0,1]$. Then put:
\begin{equation} \label{uniform1}
 k \lhd k^\prime
\qquad \text{ iff  } \qquad U_k < U_{k^\prime}.
\end{equation}
  It is immediate that the restriction $\lhd^{(n)}$
is uniform whence, by the uniqueness property of the
projective limit, $\lhd$ has $\lambda$ as its distribution.
Notice that, for each $k$, \begin{equation} \label{uniform2}
U_k= \lim_{n\rightarrow \infty} \frac{1}{n}\sum_{i\leq n, i
\neq k} 1_{(i \lhd k)} \qquad\qquad a.s., \end{equation} and
thus we may regard $U_k$ as  being the (relative) position of
the card carrying label $k$. A slightly different way of
thinking about this:  equations \eqref{uniform1} and
\eqref{uniform2} set up a measure isomorphism between the
space of total orderings endowed with $\lambda$  and the space
$[0,1]^\infty$ endowed with the infinite product of uniform
measure on $[0,1]$. It's often much easier to think about
things in the $[0,1]^\infty$ world, as we shall see.

  Suppose that $r$ is an arbitrary bijection of ${\mathbb N}$
  onto
itself and define the induced relabelling $\hat{r}$ via 
\begin{equation}
k \stackrel{\hat{r}}{\lhd} k^\prime \qquad\text{ iff } \qquad
r(k) \lhd r(k^\prime). \end{equation} It is easy to see that
$\hat{r}$ preserves the uniform measure $\lambda$ and in fact
this invariance property characterises $\lambda$.
 It is also true  that if we use \eqref{uniform2} to define a
random variable $U_k$ on the space of orderings then
\begin{equation} \label{aut} U_k\circ \hat{r}= U_{r(k)}
\end{equation} almost surely under $\lambda$. This  means that
the action of a relabelling on the space $[0,1]^\infty$ is
just to permute the co-ordinates.

To see equation \eqref{aut}, just note that given  $r$, 
for $\lambda$ almost all $\lhd$, we may define a  new order
$\lhd^\prime$  by $k\lhd^\prime k^\prime$ iff $U_{r(k)}<
U_{r(k^\prime)}$ which agrees with $\stackrel{\hat{r}}{\lhd}$.
 Some attention should be paid to the null sets here. Equation
\eqref{aut} holds except for a null set that depends on $r$.
In fact any two orderings that are dense  and open (so between
any two elements there is a third and there are no minimal or
maximal elements) have the same order-type (see Fraenkel (1976)) and
so there is some relabelling carrying one to the other.Thus
for any $\lhd$ there is a choice of $r$ so that \eqref{aut}
fails to hold at $\lhd$.

Suppose that $\hat{S}$ is a  map from \( \tilde{\mathcal O}\),
the space of orderings
of ${\mathbb N}$, into itself. When is it appropriate to call
$\hat{S}$ a shuffle? When the commutatation property
\eqref{com} holds  $\lambda$ almost surely for each relabeling
of the infinite pack as defined in the previous paragraph. In
this case there exists a unique function, $S: [0,1] \mapsto
[0,1]$, which preserves Lebesgue measure, such that, for each $k$,
\begin{equation} \label{det} U_k \circ \hat{S}= S \circ U_k
\end{equation} $\lambda$ almost surely. Moreover, each such $S$
corresponds to some $\hat{S}$.

We define a random shuffle as a suitable generalisation of such
functions:
\begin{definition}
A random shuffle is
described by a family of transition kernels $\kappa$ on the
space of orderings of ${\mathbb N}$, satisfying the following
generalisation of the commutation relation for any
$r$: whenever $A$ is a measurable subset of the space of total
orderings, and $\hat{r}$ a relabelling, \begin{equation}
\label{com2} \kappa(\hat{r}(\lhd),\hat{r}(A))= \kappa(\lhd, 
A) \qquad
 \text{ for $\lambda$ almost all } \lhd. 
\end{equation}
\end{definition}
Here we have written $\hat{r}(\lhd)$ for the ordering
$\stackrel{\hat{r}}{\lhd}$. As with deterministic shuffles, we
can express $\kappa$ using the card positions.

\begin{theorem}\label{represent}
Suppose that $\nu$ is a probability  measure on $[0,1]^2$
having both marginals uniform on $[0,1]$. Take a sequence of
independent pairs of  random variables \\
 $\bigl((U_1,V_1) \ldots , (U_k,V_k), \ldots\bigr)$, each pair
 distributed according 
to $\nu$. This then determines, by virtue of \eqref{uniform1},
the joint law of a pair of orderings $(\lhd,\lhd^\prime)$.
Take $\hat{\nu}(\lhd, \cdot)$ to be a regular conditional
probability for $\lhd^\prime$ given $\lhd$. Then
$\kappa=\hat{\nu}$  satisfies the commutation relation
\eqref{com2}. Moreover, any $\kappa$ satisfying the relation 
\eqref{com2} is a mixture of kernels constructed in this
manner. \end{theorem} \begin{proof} Suppose that
$(U_k,V_k)$ for $k \geq 1$ form a sequence of independent
pairs of random variables, each pair having the distribution
$\nu$ on $[0,1]^2$. Then the sequence of independent uniform
variables $(U_k; k \geq 1)$ gives rise to, with probability
one, an ordering $\lhd$ distributed according to $\lambda$,
and similarly  $(V_k; k \geq 1)$ gives rise to an ordering
$\lhd^\prime$. Fix a relabelling $r$. Since the ordering
$\hat{r}(\lhd)$ corresponds to the sequence of random
variables $\tilde{U}_k= U_{r(k)}$ and  similarly the ordering
$\hat{r}(\lhd^\prime)$ corresponds to the sequence of random
variables $\tilde{V}_k= V_{r(k)}$ we see that: \[
\bigl(\hat{r}(\lhd),\hat{r}(\lhd^\prime)\bigr)\stackrel{{\text
law}}{=}\bigl(\lhd,\lhd^\prime\bigr). \] From this it follows
that $\hat{\nu}$, defined as a regular conditional probability
for $\lhd^\prime$ given $\lhd$, satisfies \eqref{com2}.

 To see the last claim of the theorem, suppose that
 $\kappa$ satisfies 
\eqref{com2}, and consider a pair of orderings
$(\lhd,\lhd^\prime)$ determined as follows. Let $\lhd$ be
distributed according to $\lambda$, and the let the
conditional distribution of $\lhd^\prime$ given $\lhd$ be
$\kappa(\lhd, \cdot)$. It follows from the invariance of
$\lambda$ under relabellings and \eqref{com2} that, for any
relabelling $r$, \[
\bigl(\hat{r}(\lhd),\hat{r}(\lhd^\prime)\bigr)\stackrel{{\text
law}}{=}\bigl(\lhd,\lhd^\prime\bigr). \] Now, as we remarked
above, $\lambda$ is characterized by its invariance under
relabellings and so the law of $\lhd^\prime$ must also be
$\lambda$.
 Let the card positions corresponding to $\lhd$ be  $(U_1,
 \ldots, U_k,\ldots)$, and those corresponding to 
$\lhd^\prime$ be $(V_1, \ldots, V_k,\ldots)$.  Then the
sequence of pairs $((U_1,V_1), \ldots (U_k,V_k),\ldots)$ is
exchangeable. So the sequence admits a random empirical
measure $\Xi$ on $[0,1]^2$. Let the law of $\Xi$ be \[ \int
\alpha(d\nu) \nu, \] the integral being with respect to a
probability measure $\alpha$ on the space of probability measures
on $[0,1]^2$. Applying the strong law of large numbers to each
of the sequences $U_k$ and $V_k$, we deduce that the marginals
of $\Xi$ are, with probability one, uniform. Thus $\alpha$
must be supported on the set of  measures $\nu$ having uniform
marginals. Finally observe that, since, conditional on
$\Xi=\nu$, the pairs $(U_k,V_k)$ are independent and distributed
as $\nu$,
\[ \kappa(\lhd, \cdot)= \int \alpha(d\nu)
\hat{\nu}(\lhd, \cdot), \]
for $\lambda$ almost all $\lhd$. 
 By choosing an  appropriate version for $\hat{\nu}$ we can
 obtain equality for all $\lhd$. \end{proof}

Notice how this relates  to  the representation, \eqref{det},
of deterministic shuffles. Corresponding to  a
measure-preserving $S:[0,1] \mapsto [0,1]$ is the $\nu$ with
\begin{equation} \label{detnu} \nu(A\times B)= \int_A
1_B(S(x))dx,\hbox{ for any Borel subsets $A$ and $B$
of $[0,1]$.} \end{equation}

The above discussion seems to be the end of the story but let
us reflect.  Shuffling an infinite set of cards was really a
two step procedure. First we built the pack as the limit of a
consistent family of finite packs, then we discussed
appropriate transformations of the limiting object as
shuffles. But we could do this differently. In what follows we
consider transformations on the finite packs first, look for
some consistency of the resulting processes, and then we pass
to the limit.

Suppose that  $(\rho^{(n)}_h; h \geq 0)$ is a random walk on
${\cal S}_n$, starting from a uniformly chosen $\rho^{(n)}_0$.
Think of this as describing the state of a pack of $n$ cards
at times $h=0,1, 2, \ldots$. Now let $m<n$ and imagine that
only the cards carrying the labels $1,2 \ldots m$ are
observed. Recall that, via \eqref{correspond}, $\rho^{(n)}_h$
determines an ordering $\lhd^{(n)}_h$ and let $\lhd^{(m)}_h$
be the restriction of this ordering to $1,2 \ldots m$. Then,
using \eqref{correspond} again, we associate with
$\lhd^{(m)}_h$ a permutation  $\rho^{(m)}_h$ belonging to
${\cal S}_{m}$. Clearly, for each $h$ we have that
$\rho^{(m)}_h$  is uniformly distributed but it is easy to
construct examples so that the process $(\rho^{(m)}_h; h \geq
0)$ is not a random walk. What are the weakest conditions 
that must be placed on the jump distribution of $\rho^{(n)}$
to ensure that it {\it is} a random walk?   We do not know. But here are two
special  cases when it works.

{Case 1:}  $\lhd^{(m)}_{h+1}$ is conditionally
independent of $\lhd^{(n)}_h$ given $\lhd^{(m)}_h$.

{ Case 2:}   $\lhd^{(m)}_{h-1}$ is conditionally
independent of $\lhd^{(n)}_h$ given $\lhd^{(m)}_h$.

It is immediate that if case 1 holds then $(\rho^{(m)}_h; h
\geq 0)$ is a random walk, and, of course, case 2 is just case
1 run backwards! 

Now what we really want to do is construct an infinite family
of {\em random walks} $\bigl( (\rho^{(n)}_h; h \geq 0); n \geq
1\bigr)$ so that  the associated orderings $\lhd^{(n)}_h$ are
consistent, that is,  $\lhd^{(m)}_h$ is the restriction of
$\lhd^{(n)}_h$ whenever $m<n$.  Such a consistent family of
processes determines a limiting process $(\lhd_h; h \geq 0)$
taking values in the space of orderings of ${\mathbb N}$. 

Such a process
is Markovian with a transition  
 kernel $\kappa$ which
satisfies \eqref{com2}.
\begin{definition}
We shall call
the kernel of such a limiting process a {\it shuffle imbedding shuffle},
or SIS.
\end{definition}
When we express $\kappa$ as a
mixture:
\begin{equation}
\label{kappa} \kappa(\lhd, \cdot)=
\int \alpha(d\nu) \hat{\nu}(\lhd, \cdot),
\end{equation} the
measure $\alpha$ is suported on $\nu$ having a special form.
We investigate this special form next with the help of ${\cal I}$-invariant
orderings.

\begin{definition} Let us say
that  $(\lhd_h;h \geq 0)$ is of type 1 if case 1 holds  for
each pair $m<n$, and let us say it is of type 2 if case 2 holds
for each pair $m<n$. We always assume  that $\lhd_h$ is
distributed according to $\lambda$. We call the corresponding kernels type 1
and type 2 shuffles.

\end{definition}
\begin{remark}\label{1kernel}
It is obvious from the preceding discussion that if $\kappa$ is a type 1 shuffle then
\( \kappa(\lhd, \{ \tilde \lhd:\, \tilde\lhd^{(n)}=\sigma\}\) is
$\lambda$ almost surely constant over \( \{\lhd:\, \lhd^{(n)}=\rho\}\). We denote
the common value by $\kappa_n(\rho,\sigma)$. It is clear from the
preceding comments that $\kappa$ is determined by the $(\kappa_n)_{n
\geq 1}$.
\end{remark}
We now show that:

\begin{proposition}\label{equiv}
Each type 1 shuffle 
induces (the law of) an ${\cal I}$-invariant ordering and vice versa.
\end{proposition}

\begin{proof}Suppose that the type 1 shuffle is $\kappa$. Then we obtain the law of an
${\cal I}$-invariant ordering $\lhd$ of ${\mathbb Z}$, which we denote $P(\kappa)$, as
follows. Given integers $k_1< k_2 <\ldots <k_n$ and a
permutation $\rho \in {\cal S}_n$  then the probability that
$\lhd$ orders $k_i$ so that \begin{equation} k_i\lhd k_j
\qquad \text{ iff }\qquad \rho(i) <\rho(j) \end{equation} is
the probability that
$\rho^{(n)}_{h+1}(\rho^{(n)}_{h})^{-1}=\rho$.   Notice that
$\lhd$ is defined in such a way as to be automatically ${\cal
  I}$-invariant. In checking that this  definition  of $\lhd$
is
meaningful we need the conditional independence asserted by
case 1.

Conversely, according to  Theorem \ref{orderth}, the law, $\mathbb P$, of
$\lhd$ is a mixture: \begin{equation} \int \theta(d\mu)
{\mathbb P}^\mu, \end{equation} for some probability measure
$\theta$ on the space of probability measures on $[0,1]$. In fact,
$\theta$ is the law of the random empirical measure, $\mu(X)$, of Theorem \ref{orderth},
and is
supported on the set of quasi-uniform measures. We induce
a kernel $\Kappa({\mathbb P})$, in terms of $\theta$, as follows. Let $\mu$ be any
quasi-uniform measure  and let $X$ and $Y$ be a conjugate pair
of random variables with the law of $X$ being $\mu$. Let $U$
be an independent random variable uniformly distributed on
$[0,1]$,
  and let $\nu_\mu$
be the law of the pair $(U, UX+(1-U)Y)$. It is easy to check that
the measure $\nu_\mu$ has
uniform marginals and so, by Theorem \ref{represent}, there is a corresponding kernel,
satisfying \eqref{com2}, which we denote by $\hat{\nu}_\mu$.
Now define $\Kappa({\mathbb P})$ by
\[ \int\theta(d\mu) \hat{\nu}_\mu. \]
It remains to check that $\Kappa({\mathbb P})$ is a type 1 shuffle.
Recall that $\hat{\nu}_\mu$ is the regular conditional probability law
for $\lhd'$ given $\lhd$, where $\lhd$ and $\lhd'$ are the orders of $(U_n)_{n\geq 1}$
and $(V_n)_{n\geq 1}$ respectively and the $(U_n,V_n)_{n\geq 1}$ are
iid with common law $\nu_{\mu}$. The result now follows from the fact
(which we leave to the reader to check)
that the order of
$V_1,\ldots,V_n$ is independent of $(U_k)_{k\geq 1}$, conditional on the
order of $U_1,\ldots,U_n$.
\end{proof}

\begin{proposition}\label{equiv2} There is a one-to-one correspondence
between the laws of ${\cal I}$-invariant orderings
 and the type 1 shuffles. Under this bijection, the law of
 the ordering \[ \int \theta(d\mu) {\mathbb P}^\mu, \]
 corresponds to the kernel \[ \int
 \theta(d\mu) \hat{\nu}_\mu. \]
\end{proposition}
\begin{proof}From the proof of Proposition \ref{equiv}, all that remains
is to establish that the maps $P$ and $\Kappa$ satisfy
\begin{equation} \label{id1}
P\circ\Kappa=id
\end{equation}
and
\begin{equation} \label{id2}
\Kappa\circ P=id.
\end{equation}
To establish \eqref{id2}, given a type 1 shuffle, $\kappa$, set $\hat
\kappa=\Kappa\circ P(\kappa)$. Recall from Remark \ref{1kernel} that $\kappa$ is characterised by
$(\kappa_n(id,\rho);\rho\in S_n,n\geq 1)$ and observe that
\[\kappa_n(id,\rho)={\mathbb P}(\lhd^{(n)}=\rho)=\int \theta(d\mu){\mathbb P}^{\mu}(\lhd^{(n)}=\rho)
=\int \theta(d\mu)\hat{\nu}_\mu(id,\rho)=\hat \kappa_n (id,\rho),\]
where ${\mathbb P}=P(\kappa)$.
The proof of \eqref{id1} is similar.
\end{proof}
\begin{remark}
The proof of Proposition \ref{equiv} now makes clear the role of quasi-uniform
pairs \( \mu, \mu'\) in constructing 
type 1 shuffles. An extremal type 1 shuffle is realised by taking the
appropriate quasi-uniform \( \mu \), constructing a corresponding sequence of conjugate iid
pairs $(X_n,Y_n)$ and then setting $V_n=U_nX_n+(1-U_n)Y_n$, where $U_n$
and $V_n$ are, respectively, the initial and final positions of card
$n$. This definition still makes sense even if there are
ties in final card positions. For suppose that $V_n=V_m$; notice that this can only happen if either the
corresponding initial positions are the same and the corresponding
conjugate pairs \( (X_n,Y_n)\) and \( (X_m,Y_m)\) are equal or if the
initial positions take values in \( \{0,1\} \) and the conjugate pairs
lie on adjacent components of $G$.
In the latter case we resolve the tie by ordering $m$ above $n$ iff \(
(X_m,Y_m)\) belongs to the higher/rightmost component of $G$. In the former case,
we preserve the initial ordering between $m$ and $n$ if
$Y_n=Y_m<X_m=X_n$ and otherwise reverse it (just as in \ref{condition}).
The corresponding kernel, $\nu_{\mu}$ is defined on all $\lhd\in \tilde{\mathcal O}$, and
$\nu_{\mu}(\lhd,\cdot)$ is a probability measure on $\tilde{\mathcal O}$ for
every $\lhd$. Under ${\mathbb P}^\mu$, the law of the restriction of \( \lhd\) to \( \mathbb N\) is equal to
$\nu_{\mu}(id,\cdot)$.

\end{remark}

For type 2 shuffles the story is similar. An ${\cal 
I}$-invariant ordering is determined as follows. For integers
$k_1< k_2 <\ldots <k_n$ and a permutation $\rho \in {\cal
S}_n$  then the probability that $\lhd$ orders $k_i$ so that
\begin{equation} k_i\lhd k_j \qquad \text{ iff }\qquad \rho(i)
<\rho(j) \end{equation} is the probability that
$\rho^{(n)}_{h}(\rho^{(n)}_{h+1})^{-1}=\rho$. The law of this
ordering determines  the transition kernel as in the proof of Proposition \ref{equiv}.
 If $\mu$ is a quasi-uniform measure, let $\nu^\mu$ be the
 measure on $[0,1]^2$ defined by
\begin{equation}
\nu^\mu(dx,dy)= \nu_\mu(dy,dx),
\end{equation}
and let $\hat{\nu}^\mu$ be the associated kernel.
\begin{proposition}
 There is a one-to-one correspondence between the laws of
 ${\cal I}$-invaniant orderings and type 2 shuffles.
Under this bijection, the law of the ordering \[ \int
 \theta(d\mu) {\mathbb P}^\mu, \] corresponds to the
kernel
\[ \int \theta(d\mu) \hat{\nu}^\mu. \]
\end{proposition}
\begin{proof}The result follows immediately from Propositions
\ref{equiv} and \ref{equiv2} by time reversal.
\end{proof}

It is this result which generalizes the GSR shuffle. In particular, if the
Markov chain  corresponds to a measure $\theta$ which puts all its
mass on a single quasi-uniform  measure $\mu$, and $\mu$ is purely
atomic, then there exists a function $S:[0,1] \mapsto [0,1]$
such that \eqref{detnu} holds for $\nu^\mu$, and the chain is
in fact deterministic and is obtained by iterating the shuffle
$\hat{S}$ associated by equation \eqref{det} with $S$. If $G$
corresponding to $\mu$ is decomposed as
\[ G=\bigcup_n (l_n,L_n)\cup\bigcup_m (r_m,R_m),\]
where the atoms of $\mu$ are situated on $\{l_n;\, n\geq 1\}\cup\{R_m;\,
m\geq 1\}$ then, at least for $x\in G$,
\[ S(x)=\sum_n \frac{L_n-x}{L_n-l_n}1_{(l_n,L_n)}(x)+\sum_m \frac{x-r_n}{R_n-r_n}1_{(r_n,R_n)}(x).\]
Thus, for
example, the GSR shuffle corresponds to the quasi-uniform measure with
atoms of size $\half$ at $\half$ and 1.

We end with the following:
\begin{conjecture}
Do there exist  Markov processes, other than those constructed
above, on the space of orderings of ${\mathbb N}$ such that,
for each $n$, the restriction of the ordering to $\{1, \ldots
n\}$ evolves as if induced by a random walk on ${\cal S}(n)$?
\end{conjecture}

 \end{document}